\documentclass[preprint,12pt]{elsarticle}

\usepackage{amsmath}
\usepackage{amssymb,color}
\usepackage{subfig}
\usepackage{mathtools}
\usepackage{graphicx}
\newtheorem{Theorem}{Theorem}
\newtheorem{Lemma}{Lemma}

\newtheorem{Notation}{Notation}

\newtheorem{Remark}{Remark}
\newtheorem{ass}{Assumption}[section]

\journal{}

\begin{document}

\begin{frontmatter}

%% Title, authors and addresses
%% \title{Title\tnoteref{label1}}
%% \tnotetext[label1]{}
%% \author{Name\corref{cor1}\fnref{label2}}
%% \ead{email address}
%% \ead[url]{home page}
%% \fntext[label2]{}
%% \cortext[cor1]{}
%% \address{Address\fnref{label3}}
%% \fntext[label3]{}

\title{Superconvergence analysis of FEM and SDFEM on graded meshes for a  problem with characteristic layers %\tnoteref{t1}
	}
%\tnotetext[t1]{This paper has been supported by the Ministry of Education, Science and Technological Development of the Republic of Serbia, project no.
%	451-03-68/2020-14/200134 and project no. 451-03-68/2020-14/200156: "Innovative scientific and artistic research
%	from the FTS activity domain".}
\author[MBr]{M. Brdar\corref{cor1}}
\ead{mirjana.brdar@uns.ac.rs}
\author[GR]{G. Radojev}
\ead{radojev@dmi.uns.ac.rs}
\author[HGR]{H. -G. Roos}
\ead{hans-goerg.roos@tu-dresden.de}
\author[LjT]{Lj. Teofanov}
\ead{ljiljap@uns.ac.rs}
\address[MBr]{Faculty of Technology, University of Novi Sad,
Bulevar cara Lazara 1, 21000 Novi Sad, Serbia}
\address[GR]{Department of Mathematics and Informatics, Faculty of Sciences, University of Novi Sad, Trg Dositeja Obradovi\'{c}a 4, 21000 Novi Sad, Serbia}
\address[HGR]{Institute of Numerical Mathematics, Technical University of Dresden, Dresden D-01062,Germany}
\address[LjT]{Faculty of Technical Sciences, University of Novi Sad, Trg Dositeja Obradovi\'{c}a 6, 21000 Novi Sad, Serbia}
\cortext[cor1]{Corresponding author}

\begin{abstract}
We consider a singularly perturbed convection-diffusion with exponential and  characteristic boundary layers. The problem is numerically solved by the FEM and SDFEM method with bilinear elements on a graded mesh. For the FEM we prove almost uniform convergence and superconvergence. The use of graded mesh allows for the SDFEM to prove almost uniform esimates in the SD norm, which is not possible for Shishkin type meshes.
\end{abstract}

\begin{keyword}
%% keywords here, in the form: keyword \sep keyword
 singular perturbation \sep characteristic layers \sep finite element method \sep  streamline diffusion method\sep  graded mesh \sep superconvergence
 %% PACS codes here, in the form: \PACS code \sep code

%% MSC codes here, in the form: \MSC code \sep code
%% or \MSC[2008] code \sep code (2000 is the default)
\end{keyword}

\end{frontmatter}

%----------------------------------------------------------------------
 \section{Introduction}
%----------------------------------------------------------------------

 In this paper we consider the following convection-diffusion problem
 \begin{align}\label{1}
 Lu:=-\varepsilon\Delta u -  b u_x + c u&=f \quad \mbox{in} \quad \Omega=(0,1)\times(0,1), \nonumber\\[1ex]
 u&=0 \quad \mbox{on} \quad\partial\Omega,
 \end{align}
 with
 \begin{align*}
 b\in W^{1,\infty}(\Omega), \quad c\in L^{\infty}(\Omega), \quad b\geq\beta, \quad c \geq 0 \quad \mbox{on}\,\, \bar{\Omega},
 \end{align*}
 where $\beta$ is a positive constant and $0<\varepsilon\ll 1$ a small perturbation parameter.
 We additionally assume that
 \begin{align}\label{4}
 \displaystyle c+\frac{1}{2} b_x \geq \gamma >0,
 \quad (x,y)\in\Omega,
 \end{align}
 for some constant $\gamma,$ which will ensure the coercivity of the bilinear form associated with the differential operator $L$.

Problem \eqref{1} belongs to the class of singularly perturbed problems whose solutions are characterized by the so-called boundary layers - parts of the domain where the solution changes abruptly, i.e. where the derivatives of the solution are very large. It is known that the presence of boundary layers in solutions of singularly perturbed problems makes the application of the standard numerical procedures unstable and unsatisfactory. Therefore, the construction of at least almost $\varepsilon$-uniform numerical methods, which provide accurate approximate solution in the whole domain, is the main issue in a treatment of boundary layer problems. Problem \eqref{1} can be used in modelling of flow past a surfice and the analysis of this problem and possibilities of its numerical solutions can be of help in the numerical treatment of more complex problems.

The solution of problem (\ref{1}) is characterized by an exponential layer at $x=0$ and two parabolic layers at characteristic boundaries $y=0$ and $y=1$. Parabolic layers occur if the dominant part of the solution near some part of the boundary satisfies a partial differential equation of parabolic type. The variable in the streamline direction plays the role of the time variable in this equation. This can be observed from the asymptotic expansion.
It is known that for problems with parabolic layers does not exist a fitted scheme that converges uniformly on a uniform mesh \cite{Shishkin}. Therefore, a general strategy in the construction of a numerical method for the  singularly perturbed problem with parabolic layers is to apply some finite difference or finite element method (FEM) on a specially designed layer-adapted mesh.

In this paper we will use a graded Duran-Lombardi (DL) mesh introduced in \cite{DuranLombardi}. This mesh is defined implicitly by a recursive formula. It is a simplified version of Gartland mesh from \cite{Gratland} where the recursive formula is based on equidistribution of pointwise error, and contains exponential function in it. DL mesh has the great advantage of being simple and not requiring the a priori definition of transition points. It is robust in the sense that a mesh defined for some fixed value of the perturbation parameter can also be used for lager values of the parameter. Moreover, for a certain range of $\varepsilon$ this construction is even a better option then a mesh generated with a corresponding perturbation parameter.  In the FE analysis done so far on DL meshes their simplicity was utilized to reduce initial assumptions on solution properties.

The streamline diffusion FEM (SDFEM) which adds weighted residuals to the Galerkin FEM is one of the most frequently studied and most popular stabilized FEM. This method was proposed first in \cite{HughesBrooks} and applied to various problems.  Compared with the standard Galerkin FEM, the SDFEM provides additional control over the convective derivative in the streamline direction because of the definition of the induced streamline diffusion norm. This additional bound prevents the discrete solution from oscillating over a large part of the domain. It is well known that the SDFEM has high accuracy away from layers and good stability properties. However, in layer regions the SDFEM fails to compute accurate solutions unless layer-adapted meshes are used. There are lots of results for singularly perturbed problems concerning SDFEM on Shishkin meshes; here we refer to some  of them dealing with problem \eqref{1}: \cite{FranzKelloggStynes,FranzLinssRoos,LiuZhang1,LiuZhang2,ZhangLiu}. In the analysis of the SDFEM on Shishkin type meshes it is not possible to prove the desired estimates in the SD norm, instead only an related energy norm is used. This remedy was improved by Zhang and Liu \cite{LiuZhang}, they used modified SD norm. Franz \cite{Franz} used the same trick to modify the LPS norm for the local projection stabilization method. Surprisingly, for graded meshes estimates in the original SD norm are possible, see \cite{YinZhu} for a problem with exponential layers only. We show that this can also be realized for problems with exponential and characteristic layers.

In this paper we give for the first time a detailed superconvergence analysis of the SDFEM with bilinear elements on a DL mesh for problem (\ref{1}). This analysis provides a certain choice of the streamline diffusion parameters. The optimal choice of SD parameter on anisotropic meshes is still an open problem. For one-dimensional problems, the SD parameter can be chosen in such a way that the discrete solution is exact in mesh points or that it makes the coefficient matrix an M-matrix. For a two-dimensional problem, this approach can not be applied. Furthermore, for a problem with characteristic boundary layers it is more difficult to tune the SD parameter.

The paper is organized as follows. In Section \ref{sec:solution} we set an assumption regarding solution properties of the problem (1). In Section \ref{sec:mesh} we describe layer-adapted DL mesh and give some auxiliary estimates on the solution derivatives in $L^2$ norm. Sections \ref{sec:finite} and \ref{sec:superfem} contain proofs of convergence and superconvergence  of Galerkin FEM on a DL mesh. The formulation and basic features  of SDFEM are given in Section \ref{sec:sdfem}.  Section \ref{sec:supercsdfem} is devoted to the proof of superconvergence result of SDFEM on a DL mesh under a certain choice of the SD parameter. Finally, numerical results are presented in Section \ref{sec:numerics}, and a summary of the results is given in Section \ref{sec:summary}.

% The outline of the paper is as follows. In Section

 \begin{Notation}
 For a set $D$, a standard notation
 for Banach spaces $L_p(D)$, Sobolev spaces $W^{k,p}(D)$, $H^k(D)=W^{k,2}(D)$, norms $\|\cdot\|_{L^p(D)}$ and seminorms
 $|\cdot|_{H^k(D)}$ are used. Specially, if $p=2$ we denote the norm with $\|\cdot\|_{0,D}.$ The standard scalar product in $L^2(D)$ is marked with $(\cdot,\cdot)_D$.  Throughout the paper, we often use notation $A\lesssim B$ if a generic constant
  $C$ independent of $\varepsilon$ and mesh parameter $h$ exists,
  such that $A\leq CB$.
  %For $A\lesssim B$ and $B\lesssim A$ we write shortly $A\sim B$.
 \end{Notation}

%----------------------------------------------------------------------
\section{Solution properties}\label{sec:solution}
%----------------------------------------------------------------------

 The forthcoming error analysis is based on some a priori knowledge about the behaviour of the problem solution. Therefore, in the following assumption we give a solution decomposition and bounds of their components and derivatives. The validity of this assumption for constant coefficient problem under sufficient smoothness and compatibility conditions for function $f$ is proved in \cite{KellogStynes,Roos}. Information about the location of layers obtained from these estimates is also important for the construction of a layer-adapted mesh.

 \begin{ass} \cite{FranzLinss}\label{ass}
 The solution $u$ of problem (\ref{1}) can be decomposed as $u=v+w_1+w_2+w_{12},$ where for all $x,y\in[0,1]$  and $0\leq i+j \leq k$ we have
 \begin{alignat}{2} \label{est1}
 &\left|\partial_x^i\partial_y^j v(x,y)\right| \lesssim 1,\nonumber\\[0.5ex]
 &\left|\partial_x^i\partial_y^j w_1(x,y)\right| \lesssim\varepsilon^{-i}\exp^{-\beta x/\varepsilon},\\[0.5ex]
 &\left|\partial_x^i\partial_y^j w_2(x,y)\right| \lesssim \varepsilon^{-j/2}(\exp^{-y\delta/\sqrt{\varepsilon}}+\exp^{-(1-y)\delta/\sqrt{\varepsilon}}),\nonumber \\[0.5ex]
 &\left|\partial_x^i\partial_y^j w_{12}(x,y)\right| \lesssim\varepsilon^{-(i+j/2)}\exp^{-\beta x/\varepsilon}(\exp^{-y\delta/{\sqrt{\varepsilon}}}+\exp^{-(1-y)\delta/{\sqrt{\varepsilon}}}) \nonumber
 \end{alignat}
 for some $\delta >0.$
 For $0\leq i+j \leq k+1$ the $L^2$ bounds are
\begin{alignat}{2}  \label{est2}
 ||\partial_x^i\partial_y^j v(x,y)||_{0,\Omega} &\lesssim 1, \qquad & ||\partial_x^i\partial_y^j w_1(x,y)||_{0,\Omega} &\lesssim\varepsilon^{-i-1/2},
 \\[0.5ex]
 ||\partial_x^i\partial_y^j w_2(x,y)||_{0,\Omega} &\lesssim\varepsilon^{-j/2+1/4}, \qquad & ||\partial_x^i\partial_y^j w_{12}(x,y)||_{0,\Omega} &\lesssim \varepsilon^{-i-j/2+1/4}. \nonumber
\end{alignat}
The solution $u$ of problem (\ref{1}) satisfies
\begin{alignat}{2} \label{ass1}
\left| \partial_x^i\partial_y^j u(x,y)\right| & \lesssim \left( 1 + \varepsilon^{-i} \exp^{-\beta x/\varepsilon} + \varepsilon^{-j/2}(\exp^{-y\delta/\sqrt{\varepsilon}}+\exp^{-(1-y)\delta/\sqrt{\varepsilon}}) \right. \nonumber\\[0.5ex]
& \left. + \varepsilon^{-(i+j/2)}\exp^{-\beta x/\varepsilon}(\exp^{-y\delta/{\sqrt{\varepsilon}}}+\exp^{-(1-y)\delta/{\sqrt{\varepsilon}}})\right)
\end{alignat}
for all $(x,y)\in \Omega,$  $0\leq i+j \leq k,$ and some $\delta>0.$
\end{ass}

\begin{Remark}
	Note that the estimate \eqref{ass1} follows from \eqref{est1}. In the classical analysis from \cite{DuranLombardi,DuranLombardiPrieto} on graded meshes the authors use only estimates for the derivatives like \eqref{ass1} instead a solution decomposition. We prefer to use a little bit stronger assumption for the solution decomposition to obtain a better superconvergence result for SDFEM, see Section \ref{sec:supercsdfem}.
\end{Remark}

%----------------------------------------------------------------------
\section{The graded DL mesh}\label{sec:mesh}
%----------------------------------------------------------------------

The graded DL mesh was introduced in \cite{DuranLombardi} and here we adapt it for problem \eqref{1}.  For a given parameter $0<h<1$, the mesh in $x$-direction $\Omega^h_x$ is recursively graded with mesh points given by
 \begin{equation}\label{DLx}
 \begin{array}{lll}
 x_0=0, \\[0.5ex]
 \displaystyle x_i=i h\varepsilon, & 1\leq i\leq \lceil\frac{1}{h}\rceil\\[0.5ex]
 \displaystyle x_{i+1}= x_i+h x_i,& \lceil\frac{1}{h}\rceil \leq i\leq M_x-2, \\[0.5ex]
 x_{M_x}=1,
 \end{array}
 \end{equation}
 and in $y$-direction we obtain $\Omega^h_y$ with
 \begin{equation}\label{DLy}
 \begin{array}{lll}
 y_0=0, \\[0.5ex]
  \displaystyle y_j=j h\sqrt{\varepsilon}, & 1\leq j\leq \lceil\frac{1}{h}\rceil\\[0.5ex]
  \displaystyle y_{j+1}= y_j+h y_j,& \lceil\frac{1}{h}\rceil\leq j\leq M_y-2, \\[0.5ex]
 y_{M_y}=1/2, \\[0.5ex]
 y_{M_y+j}=1-y_{M_y-j}, & j=1,\ldots,M_y,\\[0.5ex]
 \end{array}
 \end{equation}
 where integers $M_x$ and $M_y$ are such that the following inequalities are valid
 \begin{alignat*}{4}
 x_{M_x-1} & <1 & \qquad \text {and} \qquad  & x_{M_x-1}(1+h)& \geq 1,\\[0.5ex]
 y_{M_y-1} & <\dfrac12 & \qquad \text {and} \qquad  & y_{M_y-1}(1+h)& \geq \dfrac12.
 \end{alignat*}
 We also assume that the last interval $(x_{M_x-1},1)$ is not too small compared to the previous one $(x_{M_x-2},x_{M_x-1}).$ If this is not the case the node $x_{M_x-1}$ should be eliminated. Analogously in $y$-direction. Then the DL mesh  on $\bar{\Omega}$ is given by the following tensor product
\begin{equation}
\label{mesh} \Omega^h=\Omega^h_x \times
\Omega^h_{y}.
\end{equation}
\begin{figure}[t]
\includegraphics[width=0.5\textwidth]{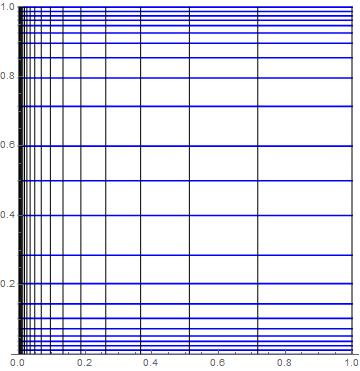}
\includegraphics[width=0.5\textwidth]{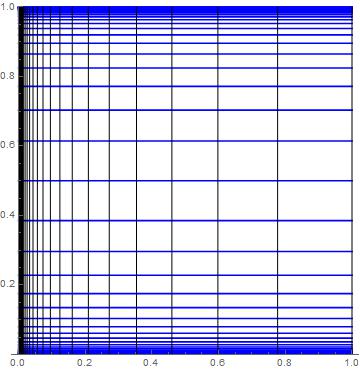}
\caption{DL mesh for $h=0.5,$
$\varepsilon=10^{-3}$ (left) and $h=0.3,$ $\varepsilon=10^{-6}$ (right).} \label{fig1}\vspace{5mm}
\end{figure}
Set $N_x=\left\lfloor\dfrac1h+1\right\rfloor$. Then $x_{N_x}=N_x h \varepsilon$,
$$x_{M_x-1}=N_x h \varepsilon (1+h)^{M_x-N_x-1}<1 \qquad N_x h \varepsilon (1+h)^{M_x-N_x}\geq 1$$
and
$$y_{M_y-1}=N_x h \sqrt{\varepsilon} (1+h)^{M_y-N_x-1}<\frac12 \qquad N_x h \sqrt{\varepsilon} (1+h)^{M_y-N_x}\geq \frac12.$$
From the above inequalities, we obtain
$$\dfrac{h}{2}\varepsilon (1+h)^{M_x-N_x-1}<\dfrac{1}{2}N_x^{-1}\leq h\sqrt{\varepsilon}(1+h)^{M_y-N_x}$$
and
$$h\sqrt{\varepsilon}(1+h)^{M_y-N_x-1}<\dfrac{1}{2}N_x^{-1}\leq \dfrac{1}{2}h\epsilon(1+h)^{M_x-N_x}.$$
A simple calculation gives a relation between $M_y$ and $M_x$, i.e.
$$\dfrac{\ln{\frac{\sqrt\varepsilon}{2}}}{\ln{(1+h)}}+M_x-1<M_y<\dfrac{\ln{\frac{\sqrt\varepsilon}{2}}}{\ln{(1+h)}}+M_x+1.$$
Let the lower and the upper bound for $M_y$ be denoted by
\begin{equation}\label{bounds_M}
 L=\left\lceil\dfrac{\ln{\frac{\sqrt\varepsilon}{2}}}{\ln{(1+h)}}\right\rceil+M_x-1, \qquad U=\left\lfloor{\dfrac{\ln{\frac{\sqrt\varepsilon}{2}}}{\ln{(1+h)}}}\right\rfloor+M_x+1.
 \end{equation}
Table~\ref{Tab:bounds} shows~\eqref{bounds_M} for fixed $\varepsilon$ and different values of $h$.
\begin{table}[!h]
\centerline{%
\begin{tabular}{c||c|c|c|c}\hline
$h$ & $M_x$ & $M_y$ & $L$& $U$  \\ \hline
0.100 & 107 & 52 & 51 & 52 \\
0.095 & 113 & 54 & 54 & 55 \\
0.090 & 136 & 57 & 56 & 57 \\
0.085 & 125 & 60 & 60 & 61 \\
0.080 & 133 & 64 & 64 & 65 \\
0.075 & 141 & 68 & 67 & 68 \\
0.070 & 151 & 73 & 72 & 73 \\
0.065 & 162 & 78 & 77 & 78 \\
0.060 & 175 & 84 & 84 & 85 \\
0.055 & 191 & 92 & 92 & 93 \\
0.050 & 209 & 125 & 124 & 125\\
\hline
\end{tabular}}
\caption{\label{Tab:bounds}%
         Lower and upper bounds $L$ and $U$ for $\varepsilon=10^{-4}$}
\end{table}
 The mesh sizes $h_{x,i}=x_i-x_{i-1}$, $i=1,\ldots,M_x$, have the following properties
 \begin{equation}
 \begin{array}{lll}\label{6}
 h_{x,i}= h\varepsilon,&  1\leq i \leq \lceil\frac{1}{h}\rceil\\[0.5ex]
 h_{x,i} \leq hx,\,\,\, x\in[x_{i-1},x_i],& \lceil\frac{1}{h}\rceil+1\leq i\leq M_x.
 \end{array}
 \end{equation}
 The mesh sizes $h_{y,i}=y_i-y_{i-1}$, $i=1,\ldots,2M_y$ satisfy
 \begin{equation}
 \begin{array}{lll}\label{7}
 h_{y,j}= h\sqrt{\varepsilon},&  1\leq j\leq \lceil\frac{1}{h}\rceil, \quad 2M_y-\lceil\frac{1}{h}\rceil \leq j \leq 2M_y\\[0.5ex]
 h_{y,j} \leq hy,\,\,\, y\in[y_{j-1},y_j],& \lceil\frac{1}{h}\rceil+1\leq j\leq 2M_y-\lceil\frac{1}{h}\rceil-1.
 \end{array}
 \end{equation}

 %%%%%%%%%%%%%%%%%%%%%%%%%%%%%%%%

For the most layer-adapted meshes constructed based on a priori given information, the number of mesh points is given in advance. For a DL mesh  there is no unique parameter $h$ that generates a mesh with a fixed number of mesh nodes $N$.
To compare numerical results on a DL mesh with  results obtained on other meshes, we construct a DL mesh which gives unique $h$, for a given $N$. One additional condition should be imposed in each direction to obtain this property. For a fixed $N$, the parameters $h_x$ and $h_y$ will be calculated such that in~\eqref{DLx} and \eqref{DLy} we have
\begin{equation}\label{DLx-m}
M_x=N, \qquad x_N=x_{N-1}+h_x x_{N-1}=1
\end{equation}
\begin{equation}\label{DLy-m}
M_y=\dfrac{N}{2}, \qquad y_{N/2}=y_{N/2-1}+h_y y_{N/2-1}=\dfrac12
\end{equation}
respectively.
%The existence of these parameters is proved in~\cite{RB2019}.

%%%%%%%%%%%%%%%%%%%%%%%%%%%%

 Since the widths of characteristic and exponential boundary layers are $\cal{O}(\sqrt{\varepsilon}|\ln \varepsilon|)$ and  $\cal{O}(\varepsilon |\ln \varepsilon|)$ respectively (\cite{Kopteva},\cite[p.274]{RoosStynesTobiska1}), the domain $\bar{\Omega}$ is consequently divided into the following subdomains:
 \begin{align*}
 \Omega_1&=\bigcup\{R_{ij}:x_{i-1}<c_1\varepsilon |\ln \varepsilon|\},\\[0.5ex]
 \Omega_2&=\bigcup\{R_{ij}:x_{i-1}\geq c_1\varepsilon|\ln \varepsilon|,\,\, y_{j-1}<c_2\sqrt{\varepsilon}|\ln \varepsilon|\,\,\mbox{or}\,\,y_{j-1}>1-c_2\sqrt{\varepsilon}|\ln \varepsilon|\},\\[0.5ex]
 \Omega_3&=\bigcup\{R_{ij}:x_{i-1}\geq c_1\varepsilon |\ln \varepsilon|,\,\,c_2\sqrt{\varepsilon}|\ln \varepsilon|\leq y_{j-1}\leq 1-c_2\sqrt{\varepsilon}|\ln \varepsilon|\},
 \end{align*}
 where $R_{ij}=[x_i,x_{i-1}]\times[y_j,y_{j-1}]$ and constants $c_1$ and $c_2$ are such that
 \begin{align*}
 \Big|\frac{\partial^{i+j}u}{\partial x^i \partial y^j}\Big|\leq C\,\,\,\mbox{for}\,\,\,0\leq i+j\leq 3,\,\,\,\\ \mbox{if}\,\,\,
x>c_1\varepsilon |\ln \varepsilon| \quad \mbox{and} \quad c_2\sqrt{\varepsilon}|\ln \varepsilon|<y<1-c_2\sqrt{\varepsilon}|\ln \varepsilon|. \nonumber
 \end{align*}
 In view of Assumption \ref{ass}, it is enough to take
\begin{equation} \label{C1C2}
\displaystyle c_1>\frac{3}{\beta} \quad \mbox{and} \quad \displaystyle c_2>\frac{3}{2\delta}.
\end{equation}

\begin{figure}[t]
\begin{center}
\unitlength0.06in
\begin{picture}(32,32)
\linethickness{0.3mm}
 \put(0,0){\line(0,1){32}}
 \put(4.1,0){\line(0,1){32}}
  \put(32,0){\line(0,1){32}}
 \multiput(0,0)(0,0){1}{\line(1,0){32}}
 \multiput(4.1,4.1)(27.5,4.5){1}{\line(1,0){28}}
 \multiput(0,32)(0,4.5){1}{\line(1,0){32}}
 \multiput(4.1,27.5)(27.5,4.5){1}{\line(1,0){28}}
 \put(0.5,15){$\Omega_{1}$}
  \put(16,15){$\Omega_3$}
 \put(16,1.3){$\Omega_2$}
 \put(16,28.7){$\Omega_2$}
\end{picture}\vspace{3mm}
\caption{Partitioning of the domain $\Omega$} \label{fig:f1}
\end{center}
\end{figure}
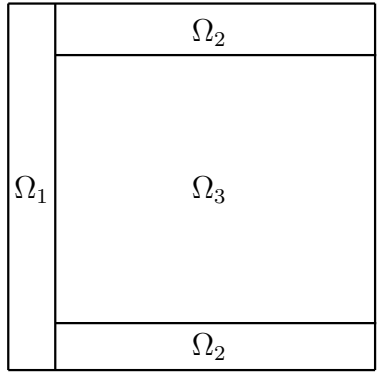

 \begin{Lemma}\cite{DuranLombardi}\label{le1}
 If $0<h<1$ is a mesh parameter and $M_x$ is the number of mesh points in $\Omega^h_x$, then
 \[h\lesssim M_x^{-1}|\ln \varepsilon|.\]
 \end{Lemma}
Analogously,  $h \lesssim M_y^{-1}|\ln \sqrt{\varepsilon}|.$ Therefore, if the total number of mesh points is denoted by $M$, then
\begin{equation}\label{h}
h\lesssim \frac{1}{\sqrt{M}}|\ln \varepsilon|.
\end{equation}

In the following two lemmas several a priori estimates for the solution of problem (\ref{1}) are given. They will be  frequently used in convergence and superconvergence analysis.
 \begin{Lemma}\label{le2}
 The solution $u$ of problem (\ref{1}) satisfy
 \begin{alignat*}{7}
 &||u_{xx}||_{0,\Omega}^2 \lesssim  \varepsilon^{-3}, & \; & ||u_{yy}||_{0,\Omega}^2 \lesssim  \varepsilon^{-3/2},
 &\; &||xu_{xx}||_{0,\Omega}^2 \lesssim \varepsilon^{-1}, &\; & ||x^2u_{xx}||_{0,\Omega}^2 \lesssim  1,\\[0.5ex]
 &||yu_{xy}||_{0,\Omega}^2 \lesssim \varepsilon^{-1}, &\; & ||u_{xy}||_{0,\Omega}^2 \lesssim  \varepsilon^{-3/2},
 &\; &||xyu_{xy}||_{0,\Omega}^2 \lesssim 1 , &\; & ||xu_{xy}||_{0,\Omega}^2 \lesssim  \varepsilon^{-1/2}, \\[0.5ex]
 &||y^2u_{yy}||_{0,\Omega}^2 \lesssim  1,& \; & ||yu_{yy}||_{0,\Omega}^2 \lesssim    \varepsilon^{-1/2}.
 \end{alignat*}
 \end{Lemma}
 %\noindent {\bf Proof}: It is easy to check that the proof follows from Assumption \ref{ass} and mesh size properties (\ref{6}) %and (\ref{7}). \qed

 \begin{Lemma}\label{le3}
 For the solution $u$ of problem (\ref{1}),  we have the following a priori estimates
 \begin{alignat*}{3}
  &\|x^3u_{xxx}\|_{0,\Omega}\lesssim 1,&\quad&\|y^3u_{yyy}\|_{0,\Omega}\lesssim 1,&\quad &\|x^2u_{xxx}\|_{0,\Omega}\lesssim \varepsilon^{-1/2},\\[0.5ex]
  &\|y^2u_{yyy}\|_{0,\Omega}\lesssim \varepsilon^{-1/4},&\quad &\|xu_{xxx}\|_{0,\Omega}\lesssim \varepsilon^{-3/2},&\quad&\|yu_{yyy}\|_{0,\Omega}\lesssim \varepsilon^{-3/4},\\[0.5ex]
 &\|u_{xxx}\|_{0,\Omega}\lesssim \varepsilon^{-5/2},&\quad&\|u_{yyy}\|_{0,\Omega}\lesssim \varepsilon^{-5/4},&\quad
 &\|u_{xxy}\|_{0,\Omega}\lesssim \varepsilon^{-7/4},\\[0.5ex]
 &\|u_{xyy}\|_{0,\Omega}\lesssim \varepsilon^{-5/4},&\quad &\|x^2u_{xxy}\|_{0,\Omega}\lesssim \varepsilon^{-1/4},&\quad&\|y^2u_{xyy}\|_{0,\Omega}\lesssim \varepsilon^{-1/2}.
  \end{alignat*}
Moreover, if (\ref{C1C2}) is satisfied then
  \begin{alignat*}{3}
   &\|u_{xxx}\|_{0,\Omega_1}\lesssim \varepsilon^{-5/2},  &\quad&\|yu_{xxy}\|_{0,\Omega_1}\lesssim \varepsilon^{-3/2},&\quad  &\|y^2u_{xyy}\|_{0,\Omega_1}\lesssim \varepsilon^{-1/2},\\[0.5ex]
    &\|x^2u_{xxx}\|_{0,\Omega_2}\lesssim \varepsilon^{1/4}|\ln\varepsilon|^{1/2},  &\quad&\|xu_{xxy}\|_{0,\Omega_2}\lesssim\varepsilon^{-1/4},&\quad  &\|u_{xyy}\|_{0,\Omega_2}\lesssim \varepsilon^{-3/4},\\[0.5ex]
 &\|x^2u_{xxx}\|_{0,\Omega_3}\lesssim 1,  &\quad&\|y^2u_{yyy}\|_{0,\Omega_3}\lesssim 1,&\quad  &\|x^2u_{xxy}\|_{0,\Omega_3}\lesssim 1,\\[0.5ex]  &\|y^2u_{xyy}\|_{0,\Omega_3}\lesssim 1,&\quad &\|xyu_{xxy}\|_{0,\Omega_3}\lesssim 1,  &\quad&\|xyu_{xyy}\|_{0,\Omega_3}\lesssim 1,\\[0.5ex]
 &\|xu_{xxx}\|_{0,\Omega_3}\lesssim 1,  &\quad&\|yu_{yyy}\|_{0,\Omega_3}\lesssim 1,&\quad
 &\|xu_{xyy}\|_{0,\Omega_3}\lesssim 1,\\[0.5ex]
 &\|yu_{xxy}\|_{0,\Omega_3}\lesssim 1.
 \end{alignat*}
 \end{Lemma}
Proofs of Lemmas \ref{le2} and \ref{le3} are based on Assumption \ref{ass} and mesh size properties (\ref{6}) and (\ref{7}). They are analogous to the  proofs of related estimates in \cite{DuranLombardi}.
%----------------------------------------------------------------------
\section{FEM} \label{sec:finite}
%----------------------------------------------------------------------
The main goal of this paper is a superconvergence result for SDFEM on a DL mesh. Therefore, the following convergence and superconvergence results for FEM on a DL mesh are necessary ingredients. The corresponding proofs are mainly based on techniques from \cite{DuranLombardi,DuranLombardiPrieto} and we present details just when they differ because of the nature of problem \eqref{1}.

 For problem \eqref{1}, the standard weak formulation is: \\find $u\in H_0^1(\Omega)$ such that $a_G(u,v)=(f,v),$ $ \forall v\in H_0^1(\Omega), $ with the bilinear form
 \begin{align*}
 a_G(w,v):=\varepsilon_1(\nabla w,\nabla v)-( b w_x+cw,v), \quad w,v\in H_0^1(\Omega).
 \end{align*}
 Let $V^h\subset H_0^1(\Omega)$ be the finite element space of
 piecewise bilinear functions defined on DL mesh \eqref{mesh}. The Galerkin finite element method is characterized by:  find $u^h\in V^h$ such that
\begin{align}\label{discret}
a_G(u^h,v^h)=(f,v^h), \quad \forall v^h\in V^h.
\end{align}
 The bilinear form satisfies the Galerkin orthogonality property. The bilinear form $a_G(\cdot,\cdot)$ is coercive with respect to the energy norm
 \begin{align}\label{energy}
 \|u\|_{\varepsilon}^2:=\varepsilon|u|_1^2+\|u\|_0^2
 \end{align}
 due to assumption \eqref{4}.
 Hence, the standard weak formulation and the Galerkin method have unique solutions.

 Given any $u\in C^0(\bar{\Omega})$ and a triangulation $T^N$ of $\Omega$ into rectangles we denote by $u^I$ the nodal piecewise bilinear interpolant to $u$ over $T^N$. Let $u^h$ be the FEM solution.

 \begin{Theorem}\label{th1}
 Let $u$ be the solution of \eqref{1}. If \eqref{ass1} with $k=2$ holds true, then on DL mesh \eqref{mesh} the interpolation error satisfies
 \begin{alignat*}{2}
 \|u-u^I\|_{0,\Omega}\lesssim h^2&, \qquad &\|u-u^I\|_{\varepsilon}\lesssim h.
 \end{alignat*}
 \end{Theorem}
 \noindent {\bf Proof}:
 We use the estimate from \cite[Theorem 3]{ApelDobrowolski} and estimates given in Lemma \ref{le2} to obtain
 \begin{alignat*}{2}
  %\|u-u^I\|_{0,R_{11}\cup R_{1,2M_y}}^2\lesssim h^4\varepsilon^{1/2}&, \qquad &
 %\sum\limits_{j=2}^{2M_y-1}\|u-u^I\|_{0,R_{1j}}^2\lesssim h^4,\\[0.5ex]
 %\sum\limits_{i=2}^{M}\|u-u^I\|_{0,R_{i1}\cup R_{i,2M_y}}^2\lesssim h^4&, \qquad & %\sum\limits_{i=2}^M\sum_{j=2}^{2M_y-1}\|u-u^I\|_{0,R_{ij}}^2\lesssim h^4,
  \sum\limits_{i=1}^{\lceil\frac{1}{h}\rceil}\sum_{j=1}^{\lceil\frac{1}{h}\rceil}\|u-u^I\|_{0,R_{ij}}^2\lesssim h^4\varepsilon^{1/2}&, \qquad &
  \sum\limits_{i=1}^{\lceil\frac{1}{h}\rceil}\sum_{j=\lceil\frac{1}{h}\rceil+1}^{M_y-\lceil\frac{1}{h}\rceil}\|u-u^I\|_{0,R_{ij}}^2\lesssim h^4&, \\[0.5ex]
 \sum\limits_{i=\lceil\frac{1}{h}\rceil+1}^{M_x}\sum_{j=1}^{\lceil\frac{1}{h}\rceil}\|u-u^I\|_{0,R_{ij}}^2\lesssim h^4&, \qquad &
 \sum\limits_{i=\lceil\frac{1}{h}\rceil+1}^{M_x}\sum_{j=\lceil\frac{1}{h}\rceil}^{M_y-\lceil\frac{1}{h}\rceil}\|u-u^I\|_{0,R_{ij}}^2\lesssim h^4&,
 \end{alignat*}
 which imply $\|u-u^I\|_{0,\Omega}\lesssim h^2.$
 Similarly, we get $\|(u-u^I)_x\|_{0,\Omega}^2 \lesssim h^2\varepsilon^{-1}$, specifically
 \begin{align} \label{intest}
 \|(u-u^I)_x\|_{0,\Omega_1}^2 \lesssim h^2\varepsilon^{-1}, \|(u-u^I)_x\|_{0,\Omega_2}^2\lesssim h^2 \sqrt{\varepsilon}|\ln \varepsilon|,  \|(u-u^I)_x\|_{0,\Omega_3}^2\lesssim h^2,
 \end{align}  and $\|(u-u^I)_y\|_{0,\Omega}^2\lesssim h^2\varepsilon^{-1/2},$ so the theorem holds true. \qed

 \begin{Theorem}\label{th2}Let $u$ be the solution of \eqref{1}. If \eqref{ass1}  with $k=2$  holds true, then for the approximate solution $u^h$ obtained by FEM with bilinear elements on DL mesh \eqref{mesh} we have
 \[
 \|u-u^h\|_{\varepsilon}\lesssim h |\ln \varepsilon| \lesssim \frac{1}{\sqrt{M}}\ln^2 \varepsilon.
 \]
 \end{Theorem}
 \noindent {\bf Proof}: Let $\eta=u-u^I,$ $\chi=u^I-u^h.$ From Galerkin orthogonality we have $a_G(u-u^h,\chi)=a_G(\eta,\chi)+a_G(\chi,\chi).$ The bilinear form $a_G(\cdot,\cdot)$ is coercive with respect to the energy norm (\ref{energy}). Let $\alpha=\min\{1,\gamma\},$ then coercivity and Galerkin orthogonality imply
  \begin{align*}
 \alpha\|\chi\|_{\varepsilon}^2\leq a_G(\chi,\chi)=-a_G(\eta,\chi)&\leq \Big| -\varepsilon(\nabla \eta, \nabla \chi)-(b\eta,\chi_x)-((c+b_x)\eta,\chi)    \Big|\\[0.5ex]
& \leq C\Big( \|u-u^I\|_{\varepsilon}\|\chi\|_{\varepsilon}+\int\limits_{\Omega}b(u-u^I)_x\chi\,d\Omega\Big)\\[0.5ex]
 & \leq C\Big(\|u-u^I\|_{\varepsilon}^2+\frac{\alpha}{2}\|\chi\|_{\varepsilon}^2+\int\limits_{\Omega}b(u-u^I)_x\chi\,d\Omega\Big).
 \end{align*}
 Then using Theorem \ref{th1}  we have
\begin{align}\label{fin1}
 \frac{\alpha}{2}\|\chi\|_{\varepsilon}^2\leq Ch^2+C\int\limits_{\Omega}b(u-u^I)_x\chi\,d\Omega
 \end{align}
 On $\Omega_1$ Poincar\'{e} inequality gives  $\displaystyle\|\chi\|_{0,\Omega_1}\leq C\varepsilon|\ln \varepsilon|\|\nabla \chi\|_{0,\Omega_1},$ using the estimates from the proof of Theorem \ref{th1}, and generalized arithmetic-geometric mean inequality we obtain
 \begin{alignat}{1}\label{fin2}
 &\int\limits_{\Omega_1}b(u-u^I)_x\chi\,d\Omega_1  \leq C\|(u-u^I)_x\|_{0,\Omega_1}\|\chi\|_{0,\Omega_1}\nonumber\\ & \leq C\|(u-u^I)_x\|_{0,\Omega_1}\varepsilon|\ln \varepsilon|\|\nabla\chi\|_{0,\Omega_1}% \nonumber\\[0.5ex]
\leq C \varepsilon\ln^2\varepsilon\|(u-u^I)_x\|_{0,\Omega_1}^2+\rho\varepsilon\|\nabla\chi\|_{0,\Omega_1}^2 \nonumber\\[0.5ex]
&\leq Ch^2\ln^2\varepsilon+\rho\varepsilon\|\nabla\chi\|_{0,\Omega_1}^2.
 \end{alignat}
 Using the same arguments on $\Omega_2$ we get
 \begin{alignat}{1}\label{fin3}
 \int\limits_{\Omega_2}b(u-u^I)_x\chi\,d\Omega_2\leq Ch^2\ln^2\varepsilon+\rho\|\nabla\chi\|_{0,\Omega_2}^2.
 \end{alignat}
 Also, on $\Omega_3$  we have
 \begin{alignat}{1}\label{fin4}
 \int\limits_{\Omega_3}b(u-u^I)_x\chi\,d\Omega_3&\leq C\|(u-u^I)_x\|_{0,\Omega_3}\|\chi\|_{0,\Omega_3}\leq C\|(u-u^I)_x\|_{0,\Omega_3}^2+\rho\|\chi\|_{0,\Omega_3}^2\nonumber\\[0.5ex]
 &\leq Ch^2+\rho\|\chi\|_{0,\Omega_3}^2.
 \end{alignat}
 For $\rho$ small enough from (\ref{fin2})-(\ref{fin4}) we obtain $\int\limits_{\Omega}b(u-u^I)_x\chi\,d\Omega \lesssim h^2\ln^2\varepsilon.$ Therefore,  (\ref{fin1}) implies  $\displaystyle\|\chi\|_{\varepsilon}\leq Ch|\ln \varepsilon|$ which together with Theorem \ref{th1} completes the proof.\qed

One should remark that Theorem \ref{th2} also holds for linear elements.

%----------------------------------------------------------------------
\section{Superconvergence of the FEM} \label{sec:superfem}
%----------------------------------------------------------------------

 In this section, we prove that the finite element approximation defined above has a property that the difference between the computed solution and the Lagrange interpolant of the exact solution is of higher order than the error itself. The forthcoming lemmas provide necessary estimates that contribute to the proof of superconvergence result.

 \begin{Lemma}\label{le4}
The solution $u^h$ of Galerkin discretization (\ref{discret}) satisfies
\begin{align*}
 \Big|\varepsilon\int\limits_{\Omega}\nabla (u-u^I)\nabla \chi dx dy\Big| \lesssim h^2\|\chi\|_{\varepsilon}\,\,\mbox{ for all}\,\, \chi\in V^h.
 \end{align*}
\end{Lemma}
 \noindent {\bf Proof}: By Lemma \ref{le3} and approach used in the proof of Theorem \ref{th1},  four estimates are obtained:
  \begin{alignat*}{1}
    \varepsilon&\sum\limits_{i=1}^{\lceil\frac{1}{h}\rceil}\sum_{j=1}^{\lceil\frac{1}{h}\rceil}
    \|(u-u^I)_xv_x\|_{0,R_{ij}}^2 \lesssim \varepsilon (h_{x,i}^2\|u_{xxx}\|_{0,R_{ij}}
    +h_{x,i}h_{y,j}\|u_{xxy}\|_{0,R_{ij}}\\[0.5ex]
 & +h_{y,j}^2\|u_{xyy}\|_{0,R_{ij}})\|v_x\|_{0,R_{ij}} \lesssim \varepsilon(h^2\varepsilon^2\|u_{xxx}\|_{0,R_{ij}}
    +h^2\varepsilon\sqrt{\varepsilon}\|u_{xxy}\|_{0,R_{ij}}\\[0.5ex] & +h^2\varepsilon\|u_{xyy}\|_{0,R_{ij}})\|v_x\|_{0,R_{ij}} \lesssim \varepsilon^{1/4} h^2\|v\|_{\varepsilon},
   \end{alignat*}
  \begin{alignat*}{1}
  \varepsilon&\sum\limits_{i=1}^{\lceil\frac{1}{h}\rceil}\sum_{j=\lceil\frac{1}{h}\rceil+1}^{M_y-\lceil\frac{1}{h}\rceil}
  \|(u-u^I)_xv_x\|_{0,R_{ij}}^2 \lesssim \varepsilon(h^2\varepsilon^2\|u_{xxx}\|_{0,R_{ij}}+h^2\varepsilon\|yu_{xxy}\|_{0,R_{ij}}\\[0.5ex]
    &+h^2\|y^2u_{xyy}\|_{0,R_{ij}})\|v_x\|_{0,R_{ij}} \lesssim  h^2\|v\|_{\varepsilon},
    \end{alignat*}
      \begin{alignat*}{1}
   \varepsilon&\sum\limits_{i=\lceil\frac{1}{h}\rceil+1}^{M_x}\sum_{j=1}^{\lceil\frac{1}{h}\rceil}
    \|(u-u^I)_xv_x\|_{0,R_{ij}}^2 \lesssim \varepsilon(h^2\|x^2u_{xxx}\|_{0,R_{ij}}+h^2\sqrt{\varepsilon}\|xu_{xxy}\|_{0,R_{ij}}\\[0.5ex]
    &+h^2\varepsilon\|u_{xyy}\|_{0,R_{ij}})\|v_x\|_{0,R_{ij}} \lesssim \varepsilon^{1/4} h^2\|v\|_{\varepsilon},
   \end{alignat*}
   \begin{alignat*}{1}
   \varepsilon&\sum\limits_{i=\lceil\frac{1}{h}\rceil+1}^{M_x}\sum_{j=\lceil\frac{1}{h}\rceil}^{M_y-\lceil\frac{1}{h}\rceil}
   \|(u-u^I)_xv_x\|_{0,R_{ij}}^2 \lesssim \varepsilon(h^2\|x^2u_{xxx}\|_{0,R_{ij}}+h^2\|xyu_{xxy}\|_{0,R_{ij}}\\[0.5ex]
    &+h^2\|y^2u_{xyy}\|_{0,R_{ij}})\|v_x\|_{0,R_{ij}} \lesssim  h^2\|v\|_{\varepsilon}.
 \end{alignat*}
  Analogously, similar estimates for $\varepsilon\|(u-u^I)_yv_y\|_{0,\Omega}^2$ are obtained, so the lemma follows.\qed

  \begin{Lemma}\label{le5} The solution $u^h$ of Galerkin discretization (\ref{discret}) satisfies
\begin{align*}
 \Big|\int\limits_{\Omega}b (u-u^I)_x \chi dx dy\Big| \lesssim h^2|\ln\varepsilon|^{\frac{1}{2}}\|\chi\|_{\varepsilon} \,\,\mbox{ for all}\,\, \chi\in V^h.
 \end{align*}
 \end{Lemma}
 \noindent {\bf Proof}: Here we use technique from \cite{FranzLinss}. In Assumption \ref{ass}  the solution decomposition $u=v+w_1+w_2+w_{12}$ is introduced. Let $\tilde{w}=w_1+w_{12},$ and
 $$\Omega_0=\bigcup\{R_{ij}:x_{i-1}<c_1\varepsilon |\ln \varepsilon|,\; \;y_{i-1}<c_2\sqrt{\varepsilon} |\ln \varepsilon| \} \subset \Omega_1.$$
Then integration by parts yields
 \[(b(u-u^I)_x,\chi)=(b(v-v^I)_x,\chi)+(b(w_2-w_2^I)_x,\chi)_{\Omega_0\cup\Omega_2}-(b_x(\tilde{w}-\tilde{w}^I),\chi)\]
 \[
 -(b(\tilde{w}-\tilde{w}^I),\chi_x)-(b_x(w_2-w_2^I),\chi)_{(\Omega_1\setminus\Omega_0)\cup\Omega_3}-
 (b(w_2-w_2^I),\chi_x)_{(\Omega_1\setminus\Omega_0)\cup\Omega_3},
 \] since $\chi\in V^h$ vanishes on the boundary of $\Omega.$  For the terms on the right-hand side we have the following estimates,
 \begin{align*}
|(b_x(\tilde{w}-\tilde{w}^I),\chi)|+|(b_x(w_2-w_2^I),\chi)_{(\Omega_1\setminus\Omega_0)\cup\Omega_3}|&\lesssim h^2\|\chi\|_{\varepsilon}\\[0.5ex]
| (b(w_2-w_2^I),\chi_x)_{(\Omega_1\setminus\Omega_0)\cup\Omega_3}|+|(b(\tilde{w}-\tilde{w}^I),\chi_x)|&\lesssim h^2\|\chi\|_{\varepsilon}\\[0.5ex]
|(b(v-v^I)_x,\chi)|&\lesssim h^2|\ln\varepsilon|^{\frac{1}{2}}\|\chi\|_{\varepsilon}\\[0.5ex]
|(b(w_2-w_2^I)_x,\chi)_{\Omega_0\cup\Omega_2}|&\lesssim h^2|\ln\varepsilon|^{\frac{1}{2}}\|\chi\|_{\varepsilon}.
 \end{align*}\qed

 \begin{Lemma}\label{le6} The solution $u^h$ of Galerkin discretization (\ref{discret}) satisfies
\begin{align*}
 \Big|\int\limits_{\Omega}c (u-u^I) \chi dx dy\Big| \lesssim h^2\|\chi\|_{\varepsilon}\,\,\mbox{ for all}\,\, \chi\in V^h.
 \end{align*}
\end{Lemma}
 The proof follows from Cauchy-Schwarz inequality and Theorem \ref{th1}.

\begin{Theorem}\label{th3}
 Let  Assumption \ref{ass} holds true. Then the FEM solution $u^h$ obtained on DL mesh \eqref{mesh} with bilinear elements and $u^I\in V^h$ satisfy
 \begin{equation} \label{sup}
 \|u^I-u^h\|_{\varepsilon}\lesssim h^2 |\ln \varepsilon|^{\frac{1}{2}}\lesssim \frac{1}{M}|\ln \varepsilon|^{\frac{5}{2}}.
 \end{equation}
 \end{Theorem}
 \noindent {\bf Proof}: The analysis starts from
 \begin{align}\label{super}
 \|u^I-u^h\|_{\varepsilon}^2\leq|a_G(\eta,\chi)|\leq\varepsilon |(\nabla\eta,\nabla\chi)|+|(b\eta_x,\chi)|+|c(\eta,\chi)|.
 \end{align}
 Above, in Lemmas \ref{le4}-\ref{le6}, we give estimates for each of the right-hand side terms. Finally, dividing (\ref{super}) by $\|u^I-u^h\|_{\varepsilon}$ the statement of theorem follows. \qed

 \begin{Remark} If we do not use the solution decomposition, but only \eqref{ass1} with $k=3$ and the technique from \cite{DuranLombardiPrieto} we obtain
 \begin{align*}
 \Big|\int\limits_{\Omega}b (u-u^I)_x \chi dx dy\Big| \lesssim h^2|\ln\varepsilon|^3\|\chi\|_{\varepsilon} \,\,\mbox{ for all}\,\, \chi\in V^h,
 \end{align*}
so the superconvergence result then is  \begin{equation} \label{sup1}\|u^I-u^h\|_{\varepsilon}\leq Ch^2|\ln\varepsilon|^3\lesssim \frac{1}{M}|\ln \varepsilon|^5.\end{equation}
 \end{Remark}

The superconvergence result can be further used to improve numerical approximation by using some postprocessing approach.

%----------------------------------------------------------------------
\section{SDFEM}\label{sec:sdfem}
%----------------------------------------------------------------------

 In order to stabilize the discretization given by the standard Galerkin FEM we introduce the streamline diffusion FEM
 \begin{align*}
 a_G(w,v)+\sum\limits_{\tau\in \Omega^N}\varrho_{\tau}(f-Lw,bv_x)_{\tau}=(f,v),
 \end{align*}
 where $\varrho_{\tau}\geq 0$ is a user chosen parameter. Its discretization reads: Find $u^h\in V^h$ such that
 \begin{align}\label{sdfem1}
 a_{SD}(u^h,v^h):=a_G(u^h,v^h)+a_{stab}(u^h,v^h)=f_{SD}(v^h)\,\,\,\mbox{for all}\,\,\, v^h\in V^h,
 \end{align}
 with
 \[a_{stab}(w,v):=\sum\limits_{\tau\in \Omega^N}\varrho_{\tau}(\varepsilon\Delta w+b w_x-cw,b v_x)_{\tau},\]
 \[f_{SD}(v):=(f,v)-\sum\limits_{\tau\in \Omega^N}\varrho_{\tau}(f,b v_x)_{\tau}.\]
 This bilinear form also satisfies the orthogonality condition.
 Now we define a streamline diffusion norm
  \[\|v\|_{SD}^2:=\|v\|_{\varepsilon}^2+\sum\limits_{\tau\in \Omega^N}\varrho_{\tau}\|b v_x\|_{0,\tau}^2.\]
 It is shown in \cite{RoosStynesTobiska1} that if
\begin{equation} \label{uslov}
0\leq\varrho_{\tau}\leq \gamma/\|c\|_{0,\tau}^2, \qquad \tau\in T^N,
\end{equation}
then $a_{SD}(v,v) \geq \frac{1}{2} \|v\|_{SD}^2.$
 Note that $\|v\|_{\varepsilon}\leq \|v\|_{SD}$ for all $v\in H_0^1(\Omega).$ Therefore, $a_{SD}(\cdot,\cdot)$ has a stronger stability then $a_{G}(\cdot,\cdot).$ Roughly,  the method is more stable when  $\varrho_{\tau}$ is closer to its upper bound.
 Problem (\ref{sdfem1}) has a unique solution $u^h\in V^h.$

 We propose the following choice for the SD parameter:
 \begin{align} \label{sdparam}
 \varrho_1\lesssim\varepsilon M^{-1}, \qquad \varrho_2\lesssim\varepsilon^{-\frac{1}{4}} M^{-1}, \qquad \varrho_3\lesssim
 \begin{dcases}
 \varepsilon^{-1}M^{-1}, & M^{-\frac{1}{2}}\leq \varepsilon\\[1ex]
 M^{-\frac{1}{2}},& M^{-\frac{1}{2}}\geq \varepsilon
 \end{dcases}
 \end{align}
 on $\Omega_1, \Omega_2$ and $\Omega_3$ respectively.

 One should remark that these bounds for the stabilization parameters are in accordance with the result of \cite{FranzLinssRoos}. Moreover, on a DL mesh it is possible to prove interpolation error estimate in the SD norm.
  \begin{Theorem}\label{thInt}
  	Let $u$ be the solution of \eqref{1}. If \eqref{ass1} with $k=2$ holds true, then on DL mesh \eqref{mesh} the interpolation error in SD norm with \eqref{sdparam} satisfies
  	\begin{alignat*}{2}
  \|u-u^I\|_{SD}\lesssim h \lesssim \frac{1}{\sqrt{M}}|\ln \varepsilon|.
  	\end{alignat*}
\end{Theorem}
The proof of Theorem \ref{thInt} follows directly from Theorem \ref{th1}, estimates \eqref{intest}, and \eqref{sdparam}.

%----------------------------------------------------------------------
\section{Superconvergence of the SDFEM}\label{sec:supercsdfem}
%----------------------------------------------------------------------
 In this section, we prove that the SDFEM approximation defined above also has a superconvergence property.
 In the proof of Theorem \ref{th4} the following lemma will be frequently used.

\begin{Lemma}\label{Prop8}\cite{FranzLinssRoos} Let $b\in W^{1,\infty}(\Omega).$ Then
	\begin{align*}
	\Big| (b(\varphi & -\varphi^I)_x,b\chi_x)_{\Omega_i}\Big| \\[0.3cm]
	&\preceq \Big[ (h_{x,i} + h_{y,i})(h_{x,i}\|\varphi_{xx}\|_{0,\Omega_i} + h_{y,i}\|\varphi_{xy}\|_{0,\Omega_i})+ h_{y,i}^2 \|\varphi_{xyy}\|_{0,\Omega_i} \Big] \|\chi_x\|_{0,\Omega_i}
	\end{align*}
	for $i=1,2,3.$
\end{Lemma}

 \begin{Theorem}\label{th4}
 Let Assumption \ref{ass} holds true. Suppose the stabilization parameter satisfies \eqref{uslov} and \eqref{sdparam}.
 Then the streamline diffusion approximation $u^h$ obtained on DL mesh with bilinear elements and $u^I \in V^h$ satisfy
 \begin{align} \label{konacna}
 \|u^I-u^h\|_{SD}\lesssim \frac{1}{M}|\ln \varepsilon|^{\frac{5}{2}}.
 \end{align}
 \end{Theorem}

 \noindent {\bf Proof}: In our error analysis we start from the coercivity and Galerkin orthogonality:
 \begin{align} \label{sd}
 \frac{1}{2}\|\chi\|_{SD}^2\leq a_G(\eta,\chi)+a_{stab}(\eta,\chi).
 \end{align}
 From the proof of Theorem \ref{th3} we have $|a_G(\eta,\chi)|\lesssim h^2|\ln \varepsilon|^{\frac{1}{2}}\|\chi\|_{\varepsilon},$
 so the second term in (\ref{sd})  has to be estimated
 \begin{align} \label{stab}
 a_{stab}(u-u^I,\chi)=\sum\limits_{\tau\in T^N}\varrho_{\tau}\Big(\varepsilon &(\Delta (u-u^I),b\chi_x)\nonumber\\[0.5ex]
 &+(b(u-u^I)_x,b\chi_x)-(c(u-u^I),b\chi_x)\Big).
 \end{align}
 We estimate these tree terms separately on different subdomains of $\Omega.$
 Estimates from Assumption \ref{ass} are frequently used in the following. For the third term in (\ref{stab}) we use Cauchy-Schwarz inequality and Theorem 1 to obtain
 \begin{align}\label{sd1}
 \varrho_2|(c(u-u^I),b\chi_x)_{\Omega_2}|&\leq C\varrho_2\|u-u^I\|_{0,\Omega_2}\|b\chi_x\|_{0,\Omega_2}\leq C\varrho_2h^2\varepsilon^{\frac{1}{4}}|\ln \varepsilon|^{\frac{1}{2}}\|b\chi_x\|_{0,\Omega_2}\nonumber\\[0.5ex]
 &\leq C\varrho_2^{\frac{1}{2}}h^2\varepsilon^{\frac{1}{4}}|\ln \varepsilon|^{\frac{1}{2}}\|\chi\|_{SD}
  \end{align}
For the second term in (\ref{stab}) let $w=w_1+w_{12}$ and $\tilde{w}=v+w_2.$
 \begin{align}\label{sd2}
 &\varrho_2|(b(w-w^I)_x,b\chi_x)_{\Omega_2}|\leq C\varrho_2\Big(\|w_x\|_{L_1(\Omega_2)}\|b\chi_x\|_{L_{\infty}(\Omega_2)}+
 \|w_x^I\|_{0,\Omega_2}\|b\chi_x\|_{0,\Omega_2}\Big)\nonumber\\[0.5ex]
 &\leq C \varrho_2\Big(\|w_x\|_{L_1(\Omega_2)}\varepsilon^{-\frac{1}{4}}|\ln \varepsilon|^{-\frac{1}{2}}\|b\chi_x\|_{0,\Omega_2}+
 \varepsilon^{\frac{1}{4}}|\ln \varepsilon|^{\frac{1}{2}}\|w_x^I\|_{L_{\infty}(\Omega_2)}\|b\chi_x\|_{0,\Omega_2}\Big)
 \nonumber\\[0.5ex]
 &\leq C\varrho_2 \Big(\varepsilon^{c_1\beta+\frac{1}{4}}|\ln \varepsilon|^{\frac{1}{2}}\|b\chi_x\|_{0,\Omega_2}+
  \varepsilon^{c_1\beta-\frac{3}{4}}|\ln \varepsilon|^{\frac{1}{2}}\|b\chi_x\|_{0,\Omega_2}\Big)\nonumber\\[0.5ex]
  &\leq C\varrho_2^{\frac{1}{2}}\varepsilon^{c_1\beta-\frac{3}{4}}|\ln \varepsilon|^{\frac{1}{2}}\|\chi\|_{SD} \leq C\varrho_2^{1/2}\varepsilon^{\frac{9}{4}}|\ln \varepsilon|^{\frac{1}{2}}\|\chi\|_{SD}
    \end{align}
based on (\ref{C1C2}).
 If we use Lema \ref{Prop8}, we get
 \begin{align}\label{sd3}
 &\varrho_2|(b(\tilde{w}-\tilde{w}^I)_x,b\chi_x)_{\Omega_2}|\nonumber\\[0.5ex]
 &\leq C\varrho_2\Big((h+h)(h\|x\tilde{w}_{xx}\|_{0,\Omega_2}+
 h\sqrt{\varepsilon}\|\tilde{w}_{xy}\|_{0,\Omega_{21}}+h\|y\tilde{w}_{xy}\|_{0,\Omega_{22}})\nonumber\\[0.5ex]
 & + h^2\varepsilon\|\tilde{w}_{xyy}\|_{0,\Omega_{21}}+h^2\|y^2\tilde{w}_{xyy}\|_{0,\Omega_{22}}\Big)\|b\chi_x\|_{0,\Omega_2}\nonumber\\[0.5ex]
  &\leq C\varrho_2^{\frac{1}{2}}\varepsilon^{\frac{1}{4}}h^2\|\chi\|_{SD},
    \end{align}
 where $\Omega_{21}$ is a part of $\Omega_2$ where $h_{y,i}$ is equidistant while $\Omega_{22}$ is a part of $\Omega_2$ where $h_{y,i}$ is non-equidistant.
 For the first part in (\ref{stab}), using H\"{o}lder inequality, we get
  \begin{align}\label{sd4}
 \varepsilon\varrho_2|(\Delta w,b\chi_x)_{\Omega_2}|&\leq C\varepsilon\varrho_2\| \Delta w\|_{L_1(\Omega_2)}\|b\chi_x\|_{L_{\infty}(\Omega_2)}\nonumber\\[0.5ex] &\leq C\varepsilon\varrho_2\varepsilon^{c_1\beta-\frac{1}{2}}|\ln \varepsilon|\|b\chi_x\|_{L_{\infty}(\Omega_2)}\nonumber\\[0.5ex]
 &\leq C\varepsilon\varrho_2\varepsilon^{c_1\beta-\frac{1}{2}}|\ln \varepsilon|\varepsilon^{-\frac{1}{4}}
 |\ln \varepsilon|^{-\frac{1}{2}}
 \|b\chi_x\|_{0,\Omega_2}\nonumber\\[0.5ex]
 &\leq C\varrho_2^{\frac{1}{2}}\varepsilon^{c_1\beta-\frac{3}{4}}|\ln \varepsilon|^{\frac{1}{2}}\|\chi\|_{SD}.
 \end{align}
 From $$(\Delta\tilde{w},b\chi_x)_{\Omega_2}+(\Delta\tilde{w},b\chi_x)_{\Omega_0}=-((b\Delta\tilde{w})_x,\chi)_{\Omega_2\cup\Omega_0},$$ where $\Omega_0=(0,\varepsilon|\ln \varepsilon|)\times(0,\sqrt{\varepsilon}|\ln \varepsilon|),$ we have
  \begin{align}\label{sd5}
 \varepsilon\varrho_2|(\Delta\tilde{w},b\chi_x)_{\Omega_2}|&\leq C\varepsilon\varrho_2\Big(\|(\Delta\tilde{w})_{x}\|_{0,\Omega_2\cup\Omega_0}\|\chi\|_{0,\Omega_2\cup\Omega_0}+
 \|\Delta\tilde{w}\|_{0,\Omega_{0}}\|\chi_x\|_{0,\Omega_0}\Big)\nonumber\\[0.5ex]
 & \leq C\varepsilon\varrho_2(\varepsilon^{-\frac{3}{4}}\|\chi\|_{0,\Omega_2\cup\Omega_0}+\varepsilon^{-\frac{1}{4}}
 |\ln \varepsilon|\|\chi_x\|_{{0,\Omega_0}})\nonumber\\[0.5ex]
  &\leq C\varrho_2\varepsilon^{\frac{1}{4}}|\ln \varepsilon|\|\chi\|_{\varepsilon}.
    \end{align}
 Collecting the above results (\ref{sd1})-(\ref{sd5}), we get
 \begin{align} \label{Om2}
 a_{stab}(u-u^I,\chi)_{\Omega_2}\leq C\Big(\varrho_2\varepsilon^{\frac{1}{4}}|\ln \varepsilon|+\varrho_2^{\frac{1}{2}}h^2\varepsilon^{\frac{1}{4}}|\ln \varepsilon|^{\frac{1}{2}}
 \Big)\|\chi\|_{SD}.
 \end{align}

For the third term in (\ref{stab}) on $\Omega_1$ we obtain
 \begin{align}\label{sd6}
 \varrho_1|(c(u-u^I),b\chi_x)_{\Omega_1}|&\leq C\varrho_1\|u-u^I\|_{0,\Omega_1}\|b\chi_x\|_{0,\Omega_1}\leq C\varrho_1^{\frac{1}{2}}h^2\|\chi\|_{SD}
  \end{align}
For the second term we proceed as follows. Let $w=w_2+w_{12}$ and $\tilde{w}=v+w_1.$
 \begin{align}\label{sd7}
 \varrho_1|(b(w-w^I)_x,b\chi_x)_{\Omega_1}|&\leq C\varrho_1\Big(\|w_x\|_{L_1(\Omega_1)}\|b\chi_x\|_{L_{\infty}(\Omega_1)}+
 \|w_x^I\|_{0,\Omega_1}\|b\chi_x\|_{0,\Omega_1}\Big)\nonumber\\[0.5ex]
   &\leq C\varrho_1^{\frac{1}{2}}\varepsilon^{-\frac{1}{2}}|\ln \varepsilon|^{\frac{1}{2}}\|\chi\|_{SD}.
    \end{align}
 Using Lemma \ref{Prop8}, the following holds
 \begin{align}\label{sd8}
 &\varrho_1|(b(\tilde{w}-\tilde{w}^I)_x,b\chi_x)_{\Omega_1}|\nonumber\\[0.5ex]
 &\leq C\varrho_1\Big(2h(h\varepsilon\|\tilde{w}_{xx}\|_{0,\Omega_{11}}+
 h\|x\tilde{w}_{xx}\|_{0,\Omega_{12}}+h\sqrt{\varepsilon}\|\tilde{w}_{xy}\|_{0,\Omega_{21}}+
 h\|y\tilde{w}_{xy}\|_{0,\Omega_{22}})\nonumber\\[0.5ex]
 & + h^2\varepsilon\|\tilde{w}_{xyy}\|_{0,\Omega_{21}}+h^2\|y^2\tilde{w}_{xyy}\|_{0,\Omega_{22}}\Big)\|b\chi_x\|_{0,\Omega_1}\nonumber\\[0.5ex]
  &\leq C\varrho_1^{\frac{1}{2}}\varepsilon^{-\frac{1}{2}}h^2\|\chi\|_{SD},
    \end{align}
 where $\Omega_{21}$ is a part of $\Omega_1$ where $h_{y,i}$ is equidistant while $\Omega_{22}$ is a part of $\Omega_1$ where $h_{y,i}$ is non-equidistant. Analogously, $\Omega_{11}$ is a part of $\Omega_1$ where $h_{x,i}$ is equidistant while $\Omega_{12}$ is a part of $\Omega_1$ where $h_{x,i}$ is non-equidistant.
 For the first part in (\ref{stab}) we get
  \begin{align}\label{sd9}
 \varepsilon\varrho_1|(\Delta w,b\chi_x)_{\Omega_1}|&\leq C\varepsilon\varrho_1\| \Delta w\|_{L_1(\Omega_1)}\|b\chi_x\|_{L_{\infty}(\Omega_1)}\leq C\varrho_1^{\frac{1}{2}}|\ln \varepsilon|^{\frac{1}{2}}\|\chi\|_{SD}
 \end{align}
and
  \begin{align}\label{sd10}
 \varepsilon\varrho_1|(\Delta\tilde{w},b\chi_x)_{\Omega_1}|&\leq C\varepsilon\varrho_1\|\Delta\tilde{w}\|_{0,\Omega_1}\|\chi_x\|_{0,\Omega_1}\leq C\varrho_1\varepsilon^{-1}\|\chi\|_{\varepsilon}.
    \end{align}
 From (\ref{sd6})-(\ref{sd10}), we get
 \begin{align} \label{Om1}
 a_{stab}(u-u^I,\chi)_{\Omega_1}\leq C\Big(\varrho_1\varepsilon^{-1}+
 \varrho_1^{\frac{1}{2}}\varepsilon^{-\frac{1}{2}}h^2|\ln \varepsilon|^{\frac{1}{2}}\Big)\|\chi\|_{SD}.
 \end{align}
On $\Omega_3$ let $w=w_1+w_2+w_{12}.$ Then Cauchy-Schwarz inequality implies
 \begin{align}\label{sd11}
 \varepsilon\varrho_3|(\Delta w,b\chi_x)_{\Omega_3}|&\leq C\varepsilon\varrho_3\| \Delta w\|_{0,\Omega_3}\|b\chi_x\|_{0,\Omega_3}\leq C\varrho_3\varepsilon^{\frac{5}{4}}\|\chi\|_{\varepsilon}
 \end{align}
 and \begin{align}\label{sd12}
 \varrho_3|(c(w-w^I),b\chi_x)_{\Omega_3}|&\leq C\varrho_3\|w-w^I\|_{0,\Omega_3}\|b\chi_x\|_{0,\Omega_3}\leq C\varrho_3^{\frac{1}{2}}h^2\|\chi\|_{SD}
  \end{align}
  The following identity from \cite{FranzLinss}
 \begin{align}\label{Lin}
 |((v-v^I)_x,\chi_x)_{\tau}|\leq C h_{y,\tau}^2\|v_{xyy}\|_{0,\tau}\|\chi_x\|_{0,\tau}\,\, {\mbox for \; all} \,\,v\in
 C^3(\bar{\tau})
 \end{align}
  applied to the second term gives
 \begin{align}\label{sd13}
 \varrho_3|(b(w-w^I)_x,b\chi_x)_{\Omega_3}|&\leq C\varrho_3h^2\|y^2w_{xyy}\|_{0,\Omega_3}\|\chi_x\|_{0,\Omega_3}
 \nonumber\\[0.5ex]
   &\leq C\varrho_3^{\frac{1}{2}}\varepsilon^{\frac{7}{4}}|\ln\varepsilon|^2 h^2\|\chi\|_{SD}\leq
   C\varrho_3^{\frac{1}{2}} h^2\|\chi\|_{SD}.
    \end{align}

    The terms containing $v$ are not exponentially
small away from the layers, so they require more careful treatment. We use
\[(\Delta v,b\chi_x)_{\Omega_3}+(\Delta v,b\chi_x)_{\Omega_1\setminus\Omega_0}=-((b\Delta v)_x,\chi)_{\Omega_3\cup(\Omega_1\setminus\Omega_0)}\] and obtain

\begin{align}\label{sd14}
 \varepsilon\varrho_3|(\Delta v,b\chi_x)_{\Omega_3}|&\leq C\varepsilon\varrho_3(\|\chi\|_{0,\Omega_1\cup\Omega_3}+\|\chi_x\|_{L_1(\Omega_1\setminus\Omega_0)})\nonumber\\[0.5ex]
 &\leq C\varepsilon\varrho_3(\|\chi\|_{\varepsilon}+(\mbox{meas}(\Omega_1))^{\frac{1}{2}}\|\chi_x\|_{0,\Omega_1\setminus\Omega_0}) \nonumber\\[0.5ex]
 &\leq C\varrho_3\varepsilon|\ln \varepsilon|^{\frac{1}{2}}\|\chi\|_{\varepsilon}.
 \end{align}

Here we also use (\ref{Lin})
\begin{align}\label{sd15}
\varrho_3|(b(v-v^I)_x,b\chi_x)_{\Omega_3}|\leq C\varrho_3 h^2\|b\chi_x\|_{0,\Omega_3}\leq C\varrho_3^{\frac{1}{2}} h^2\|\chi\|_{SD}.
 \end{align}

\begin{align}\label{sd16}
 \varrho_3|(c(v-v^I),b\chi_x)_{\Omega_3}|&\leq C\varrho_3\|v-v^I\|_{0,\Omega_3}\|b\chi_x\|_{0,\Omega_3}\leq C\varrho_3^{\frac{1}{2}}h^2\|\chi\|_{SD}
  \end{align}

Collecting (\ref{sd11})-(\ref{sd16}), we obtain
 \begin{align} \label{Om3}
 a_{stab}(u-u^I,\chi)_{\Omega_3}\leq C\Big(\varrho_3\varepsilon|\ln \varepsilon|^{\frac{1}{2}}+
 \varrho_3^{\frac{1}{2}}h^2\Big)\|\chi\|_{SD}.
 \end{align}

Applying estimates \eqref{sdparam} of the SD parameters in \eqref{Om2}, \eqref{Om1},  and \eqref{Om3}, together with \eqref{h}, we obtain \eqref{konacna}.
  \qed

%----------------------------------------------------------------------
\section{Numerical experiments}\label{sec:numerics}
%----------------------------------------------------------------------
In this section, we present some numerical experiments in order to test our theoretical results. We consider model problem from \cite{FranzLinssRoos} given by
 \begin{align*}
 -\varepsilon\Delta u-(2-x)u_x+\frac{3}{2}u&=f(x,y) \quad x\in\Omega, \\[1ex]
 u&=0, \quad \mbox{on} \quad\partial\Omega,
 \end{align*}
 where the function $f$ is chosen in such a way that
 \begin{align*}
 u(x,y)=\Big(\cos{\frac{\pi x}{2}}-\frac{e^{-x/\varepsilon}-e^{-1/\varepsilon}}{1-e^{-1/\varepsilon}}\Big)
 \frac{(1-e^{-y/\sqrt{\varepsilon}})
 (1-e^{-(1-y)/\sqrt{\varepsilon}})}{1-e^{-1/\sqrt{\varepsilon}}}
 \end{align*}
 is the exact solution.
 The rate of convergence is calculated in the standard way.
 All computations were carried out using MATLAB R2020a. We employ SDFEM with $\rho_{\tau}$ chosen to be the maximal value allowed by Theorem \ref{th4}.
\begin{figure}[t]
\begin{centering}
\includegraphics[width=0.7\textwidth]{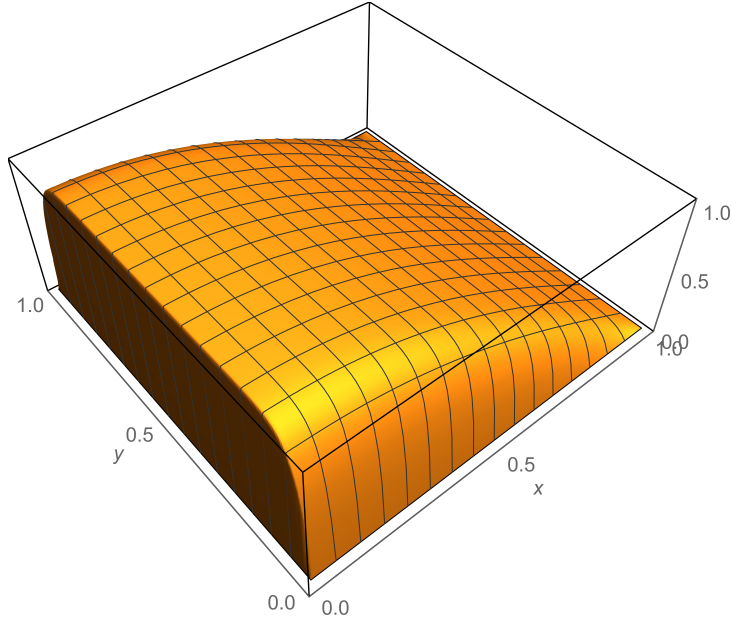}
\caption{Test problem with $\varepsilon=10^{-3}$.} \label{fig0}\vspace{5mm}
\end{centering}
\end{figure}

The errors in various norms for SDFEM are presented in Tables~\ref{tab:h}-\ref{tab:eps}.
Table~\ref{tab:h} shows the results on a DL-mesh when parameter $h$ is chosen a priori. The superconvergence property can be very well observed.
In Table~\ref{tab:1} parameter $\varepsilon$ is fixed while the number of mesh points in both directions is given in advance. Here we can clearly see the first order of convergence and the second order of superconvergence result.

In Table~\ref{tab:eps} the number of mesh points in each direction is fixed until $\varepsilon$ varies. This table shows that SDFEM on a DL-mesh is almost uniform in the singular perturbation parameter $\varepsilon$ (up to the logarithmic factor).

The last column in Tables \ref{tab:h}-\ref{tab:eps} contain the errors in $L^{\infty}$ norm  which suggest second order convergence. These results are given for the purpose of comparison with  finite difference methods. We do not have theoretical justification for it.

Finally, in Table~\ref{tab:compare} the comparison between superconvergence property of SDFEM method ona a DL mesh  and a Shishkin mesh for fixed $\varepsilon$ is given. Numerical experiments shows that SDFEM method gives better results on a DL mesh than on a Shishkin mesh for the number of mesh points greater then $32^2$.

\begin{table}[htb]
	\centerline{
		\begin{tabular}{|c|c|c|c|c|c|}  \hline
			$h$ & $M_x$  & $2M_y$ & $\|u-u^h\|_{SD}$ &     $ \|u^I-u^h\|_{SD}$ &   $\|u-u^h\|_{\infty}$  \\ \hline
			$0.15$    & 106 & 104 &   3.312e-02 &  1.728e-03 &   2.693e-03  \\ \hline
			$0.09$    & 172 & 168   &    1.997e-02  &  6.117e-04 &  1.077e-03 \\ \hline
			$0.06$    & 254 & 248 &  1.335e-02  &   2.695e-04 &   5.065e-04   \\ \hline
			$0.04$    & 378 & 368 &  8.914e-03  &  1.183e-04 &     2.334e-04    \\ \hline
			$0.02$    & 748 & 728 &  4.465e-03  & 2.971e-05 &   6.066e-05   \\ \hline
			$0.085$   & 1750 & 1704 & 1.900e-03   &    5.328e-06 &  1.113e-05 \\ \hline
	\end{tabular}}
	\caption{Errors on mesh DL~\eqref{DLx}-\eqref{DLy}  with $\varepsilon=10^{-6}$ and different $h$.}
	\label{tab:h}
\end{table}

\begin{table}[htb]
	\centerline{
		\begin{tabular}{|c||c|c||c|c||c|c|}  \hline
			{$M_x=2M_y$}  &   $\|u-u^h\|_{SD}$&  rate &  $\|u^I-u^h\|_{SD}$ &  rate &  $\|u-u^h\|_{\infty}$& rate  \\ \hline
			$16$    &   4.252e-01 &  1.31 &           2.371e-01    &  2.58 &  3.779e-01               &  2.81 \\
			$32$    &   1.715e-01 &  1.17 &           3.957e-02    &  2.25 &  5.379e-02               &  1.65 \\
			$64$    &   7.629e-02 &  1.09 &           8.334e-03    &  2.16 &  1.711e-02               &  1.77 \\
			$128$   &   3.595e-02 &  1.04 &           1.858e-03    &  2.03 &  5.012e-03               &  1.87 \\
			$256$   &   1.753e-02 &  1.02 &           4.539e-04    &  2.03 &  1.367e-03               &  1.93\\
			$512$   &   8.619e-03 &  1.01 &           1.115e-04    &  2.05 &  3.586e-04               &  1.97\\
			$1024$  &   4.276e-03 &  1.01 &           2.700e-05    &  2.01 &  9.181e-05               &  1.98\\
			$2048$  &   2.131e-03 &   -   &           6.691e-06    &   -   &  2.323e-05               &   -\\
			\hline
	\end{tabular}}
	\caption{Errors on DL mesh~\eqref{DLx-m}-\eqref{DLy-m} with $\varepsilon=10^{-8}$.}
	\label{tab:1}
\end{table}

\begin{table}[htb]
	\centerline{
		\begin{tabular}{|c|c|c||c|}  \hline
			{$\varepsilon$}
			&  $\|u-u^h\|_{SD}$ &   $ \|u^I-u^h\|_{SD}$ &   $\|u-u^h\|_{\infty}$ \\ \hline
			$10^{-2}$    &   5.456e-04  & 5.495e-06 &         4.339e-06   \\
			$10^{-3}$    &   8.370e-04  & 5.172e-06 &         1.351e-05     \\
			$10^{-4}$    &   1.107e-03  & 2.071e-06 &         4.973e-06    \\
			$10^{-5}$    &   1.367e-03  & 2.785e-06 &         4.608e-06    \\
			$10^{-6}$    &   1.623e-03  & 3.893e-06 &         7.875e-06 \\
			$10^{-7}$    &   1.877e-03  & 5.197e-06 &         1.365e-05 \\
			$10^{-8}$    &   2.131e-03  & 6.691e-06 &		  2.323e-05 \\
			$10^{-9}$    &   2.385e-03  & 8.374e-06 &		  4.012e-05 \\
			$10^{-10}$   &   2.639e-03  & 1.025e-05 &		  6.999e-05 \\
			\hline
	\end{tabular}}
	
	\caption{Errors on DL mesh ~\eqref{DLx-m}-\eqref{DLy-m} with $M_x=2M_y=2048$.}
	\label{tab:eps}
\end{table}

\begin{table}[htb]
	\centerline{
		\begin{tabular}{|c||c|c||c|c||c|c|} \hline
			& \multicolumn{2}{c||}{DL mesh}
			& \multicolumn{2}{c||}{Shishkin mesh}
			\\ \hline
			{$M_x=2M_y$}  & $ \|u^I-u^h\|_{SD}$ &  rate & $ \|u^I-u^h\|_{SD}$ &  rate   \\ \hline
			$16$    &            2.371e-01    &  2.34  &   4.857e-02 &  1.31  \\
			$32$     &           2.488e-02    &  2.25 &   1.962e-02 &  1.45 \\
			$64$    &           5.243e-03    &  2.20   &   7.170e-03 &  1.55\\
			$128$    &           1.137e-03    &  2.13 &   2.453e-03 &  1.61 \\
			$256$   &           2.600e-04    &  2.04  &   8.027e-04 &  1.66\\
			$512$   &           6.332e-05    &  2.01  &   2.542e-04 &  1.70\\
			$1024$  &           1.568e-05    &  2.01  &   7.846e-05 &  1.72\\
			$2048$   &           3.893e-06    &   -    &   2.374e-05 &   - \\
			\hline
	\end{tabular}}
	\caption{  Comparison: DL~\eqref{DLx-m}-\eqref{DLy-m} and Shishkin mesh for $\varepsilon=10^{-6}$.}
	\label{tab:compare}
\end{table}

\section{Summary} \label{sec:summary}
A singularly perturbed elliptic problem with characteristic layer has been considered.
To obtain numerical approximation of the problem, we apply Galerkin FEM and SDFEM with bilinear elements on a layer-adapted DL mesh. We  proved that such discretizations exhibit superconvergence property with the appropriate choice of streamline diffusion parameters. The results of numerical experiments confirm our theoretical results. Moreover, they show that despite its almost uniform convergence, a DL mesh can be fairly good alternative to the widely used  Shishkin mesh.

\vspace{1cm}
{\bf Acknowledgement}. This paper has been supported by the Ministry of Education, Science and Technological Development of the Republic of Serbia, project no.
	451-03-68/2020-14/200134 and project no. 451-03-68/2020-14/200156: "Innovative scientific and artistic research
	from the FTS activity domain". The authors are grateful to Professor S. Franz (TU Dresden, Germany) for his MATLAB codes.

\end{document}